\documentclass[10pt]{article}

\usepackage{amstex}
\usepackage{amssymb}
\usepackage{amscd}
\usepackage{pb-diagram}

\def\oisk{\buildrel\textstyle{\rm k \atop ^\vee} \over\dots}

\def\Talpha#1{\vbox{\ialign{##\crcr
$\alpha$\crcr\noalign{\kern2pt\nointerlineskip}
$\hfil\displaystyle{#1}\hfil$\crcr}}}

\def\Onabla#1{\vbox{\ialign{##\crcr$\,\scriptstyle{0}
$\crcr\noalign{\kern2pt\nointerlineskip}
$\hfil\displaystyle{#1}\hfil$\crcr}}}
\def\im{\mbox{Im}}

 \def\cala{{\cal A}}
\def\fraca{{\mathfrak{A}}}
\def\fracs{{\mathfrak{S}}}
\def\fracf{{\mathfrak{f}}}
\def\frach{{\mathfrak{h}}}
\def\fracm{{\mathfrak{m}}}

\def\calm{{\cal M}}

\def\calh{{\cal H}}
\def\fm{\mathfrak m}
\def\ft{\mathfrak T}
 
 \def\cals{{\cal S}}

\def\cals{{\cal S}}

\def\bbbone{\mbox{\rm 1\hspace {-.6em} l}}

\def\oiseau{\displaystyle{\buildrel \textstyle{\rm k \atop ^\vee} \over .}}

 \newtheorem{theorem}{THEOREM}
\newtheorem{lemma}{LEMMA}
\newtheorem{definition}{DEFINITION}
 
\begin{document}
\textheight=41pc
\textwidth=28pc

\begin{center} {\large\bf GENERALIZED DIFFERENTIAL SPACES WITH $d^N=0$ AND THE
$q$-DIFFERENTIAL CALCULUS} \end{center} \vspace{0.2cm}
\begin{center} Michel DUBOIS-VIOLETTE\\  \vspace{0.1cm} {\small Laboratoire de
Physique Th\'eorique et Hautes Energies\footnote{Laboratoire associ\'e au
Centre
National de la Recherche Scientifique - URA D0063}\\ Universit\'e Paris XI,
B\^atiment 211 - 91405 Orsay Cedex, France\\ flad$@$qcd.th.u-psud.fr}
\end{center}
\begin{flushright} \today \end{flushright}

{\small We present some results concerning the generalized homologies
associated with nilpotent endomorphisms $d$ such that $d^N=0$ for some integer
$N\geq 2$. We then introduce the notion of graded $q$-differential algebra and
describe some examples. In particular we construct the $q$-analog of the
simplicial differential on forms, the $q$-analog of the Hochschild differential
and the $q$-analog of the universal differential envelope of an associative
unital algebra.}

\section{Introduction}
Our aim is to discuss properties  of nilpotent endomorphisms $d$ such that
$d^N=0$ for some integer $N\geq 2$, of the corresponding generalized homologies
and of the associated generalizations of graded differential algebras. The
natural setting would be to use a category of modules carrying representations
of the group of $N$-th roots of the unity. Here however, for simplicity, we
shall work with complex vector spaces and the natural representations of this
group in such spaces, (i.e. multiplication by the corresponding complex
numbers). Such a representation is characterized by a primitive $N$-root $q$ of
the unity. For eventual applications to $q$-deformations, etc., we drop the
assumption that $q$ is a root of unity and develop the notion of
$q$-differential calculus for $q\in \mathbb C$ with $q\not=0$ and $q\not= 1$
(i.e. $q\in \mathbb C\backslash\{0,1\}$). The terminology adopted here is
influenced by a paper by M. Kapranov on similar topics \cite{Kapr}. Other
examples of spaces with $d^N=0$, etc. can be found in R. Kerner's contribution
to this volume.

\section{Generalized homology associated with $d^N=0$}

Let $E$ be a complex vector space equipped with a nilpotent endomorphism $d$
satisfying $d^N=0$ where $N$ is an integer with $N\geq 2$. One has
$\im(d^{N-k})\subset \ker(d^k)$ for $k\in \{0,1,\dots,N\}$ and therefore the
vector spaces $H^{(k)}=H^{(k)}(E)=\ker(d^k)/\im(d^{N-k})$ are well defined. In
fact $H^{(0)}=H^{(N)}=\{0\}$ and the $H^{(k)}$ are the generalization of the
homology of $E$ for $1\leq k\leq N-1$.\\

Let $\ell$ and $m$ be two positive integers such that $\ell+m\leq N$. One has
$\ker(d^m)\subset \ker(d^{\ell+m})$ and $\im(d^{N-m})\subset \im
(d^{N-(\ell+m)})$ so the inclusion $i^\ell:\ker(d^m)\rightarrow
\ker(d^{\ell+m})$ induces a homomorphism $[i^\ell]:H^{(m)}\rightarrow
H^{(\ell+m)}$. On the other hand, one has $d^m(\ker(d^{\ell+m}))\subset
\ker(d^\ell)$ and $d^m(\im(d^{N-(\ell+m)}))\subset \im(d^{N-\ell})$ and
therefore $d^m$ induces a homomorphism $[d^m]:H^{(\ell+m)}\rightarrow
H^{(\ell)}$. Notice that $[i^\ell]=[i]^\ell$ and that $[d^m]=[d]^m$.  One has
the following lemma.
\begin{lemma}
The hexagon $(\calh^{\ell,m})$ of homomorphisms
\[ \begin{diagram} \node{}
\node{H^{(\ell+m)}} \arrow{e,t}{[d^m]} \node{H^{(\ell)}}
\arrow{se,t}{[i^{N-(\ell+m)}]} \node{} \\ \node{H^{(m)}} \arrow{ne,t}{[i^\ell]}
\node{} \node{} \node{H^{(N-m)}} \arrow{sw,b}{[d^\ell]} \\ \node{}
\node[1]{H^{(N-\ell)}} \arrow{nw,b}{[d^{N-(\ell+m)}]} \node{H^{(N-(\ell+m))}}
\arrow{w,b}{[i^m]} \node{} \end{diagram} \]
is exact.
\end{lemma}
\noindent
For the proof and for more details, we refer to \cite{D-V K} and \cite{M.D-V}.
We shall use  the following criterium ensuring the vanishing of the $H^{(k)}$.

\begin{lemma}
Let $q$ be a primitive $N$-th root of the unity and assume that there is an
endomorphism $h$ of $E$ such that $h\circ d-q d\circ h=I$ where $I$ denotes the
identity mapping of $E$ onto itself. Then one has $H^{(k)}=0$ for $k\in
\{1,\dots,N-1\}$.
\end{lemma}
In fact $hd-qdh=I$ with $q$ being a primitive $N$-th root of the unity implies,
\cite{M.D-V}, $\sum^{N-1}_{k=0} d^{N-1-k} h^{N-1}d^k=[(N-1)!]_qI$, where
$[(N-1)!]_q=[(N-1)]_q \dots [2]_q$ with $[n]_q=1+q+\dots + q^{n-1}$. The result
follows since $[(N-1)!]_q \not=0$.

\subsection*{Example 1: Complex matrix algebras}

Let $N$ be as above an integer with $N\geq 2$ and let $n_k$,
$k\in\{1,\dots,N\}$ be $N$ integers greater or equal to 1, $n_k\geq 1$, with
sum $S=\sum^N_{k=1}n_k\geq N$. The algebra of complex $S\times S$ matrices will
be denoted by $M_S(\mathbb C)$.
Associated to the family $(n_k)$, there is a decomposition into rectangular
blocks $A^i_j$ of each element $A$ of $M_S(\mathbb C): A=(A^i_j)$,
$(i,j=1,2,\dots,N)$, where $A^i_j$ is a complex matrix with $n_i$ lines and
$n_j$ columns. One equips $M_S(\mathbb C)$ with a $\mathbb Z_N$-graduation,
$M_S(\mathbb C)=\oplus_{p\in \mathbb Z_N}(M_S(\mathbb C))^p$, by giving the
degree $j-i \mod(N)$ to the block $A^i_j$, (i.e. to the matrix which has only
$A^i_j$ as nonzero block). Equipped with this graduation, $M_S(\mathbb C)$ is a
$\mathbb Z_N$-graded algebra: $(M_S(\mathbb C))^a(M_S(\mathbb C))^b\subset
(M_S(\mathbb C))^{a+b}, \forall a,b\in \mathbb Z_N$. Let $e\in (M_S(\mathbb
C))^1$ be such that $e^N$ is a multiple of the unit $\bbbone \in (M_S(\mathbb
C))^0$:

\begin{equation}
e^N=\lambda \bbbone,\ \ \lambda \in \mathbb C
\label{eq1}
\end{equation}

Let $q$ be a primitive $N$-th root of the unity and set $d(A)=e A-q^aAe$ for
$A\in (M_S(\mathbb C))^a$. This defines an endomorphism $d$ of $M_S(\mathbb C)$
which is of degree 1 and satisfies
\begin{equation}
d(AB)=d(A)B+q^a Ad(B),\ \forall A\in (M_S(\mathbb C))^a,\ \forall B\in
M_S(\mathbb C)
\label{eq2}
\end{equation}
and $d^N=0$. Since $d^N=0$, the $H^{(k)}$, $k\in\{ 1,\dots, N-1\}$, are well
defined.\\
In the case $S=N$, i.e. $n_k=1$ $\forall k\in\{1,\dots,N\}$, any $e\in
(M_N(\mathbb C))^0$ satisfies (\ref{eq1}). In this case $(S=N)$, if $e^N\not=
0$ then one has $H^{(k)}=0, \forall k\in \{1,\dots, N-1\}$. $\square$

\subsection*{Example 2 : (Co)Presimplicial vector spaces}

Although all the constructions described here work as well for presimplicial
vector spaces as for copresimplicial vector spaces (dual version), we shall
describe them in the copresimplicial setting since most examples we shall meet
in the following sections are of this type. Also, it should be clear that we
work here with complex vector spaces for simplicity but that we could as well
work with an appropriate category of modules carrying representations of the
group of $N$-th roots of the unity. For the notion of (co)-presimplicial
module, we refer to J.L. Loday's book \cite{Lo}. Recall that a {\sl
copresimplicial vector space} $E$ is a family $(E^n)_{n\in \mathbb N}$ of
(complex) vector spaces together with linear mappings, {\sl the cofaces},
$\mathfrak f_k:E^n\rightarrow E^{n+1}$, $k\in \{0,\dots, n+1\}$, such that
\begin{equation}
\mathfrak f_\ell \mathfrak f_k=\mathfrak f_k\mathfrak f_{\ell-1},\ \
\mbox{for}\ \ k<\ell
\label{eq3}
\end{equation}
As usual, the dependence on $n$ of the $\mathfrak f_k$ has been dropped for
notational simplicity. We shall identify $E$ with the $\mathbb N$-graded vector
space $\oplus_n E^n$. Let $q$ be a complex number with $q\not= 0$ and $q\not=1$
and let us define two homogeneous endomorphisms of degree 1 of $E$, $d_q$ and
$\tilde d_q$, by setting for $x\in E^n$
\begin{equation}
d_q(x)=\sum^n_{k=0} q^k{\mathfrak f}_k(x)-q^n{\mathfrak f}_{n+1}(x)\
\mbox{and}\ \tilde d_q(x)=\sum^{n+1}_{k=0}q^k {\mathfrak f}_k(x)
\label{eq4}
\end{equation}
\begin{lemma}
Let $N$ be an integer with $N\geq 2$ and assume that $q$ is a primitive $N$-th
root of the unity, then one has $(d_q)^N=0$ and $(\tilde d_q)^N$.
\end{lemma}
Notice that for $N=2$, i.e. for $q=-1$, these two endomorphisms coincide with
the usual differential map of $E$, \cite{Lo}. Therefore when $q$ is a primitive
$N$-th root of the unity, the corresponding $H^{(k)}$ are generically
nontrivial. $\square$\\

In the previous example, in the case where $q$ is a primitive $N$-th root of
the unity, $E$ is $\mathbb N$-graded and $d$ is homogeneous of degree one with
$d^N=0$. More generally, let $E$ be a $\mathbb Z$-graded vector space
$E=\oplus_n E^n$ equipped with an endomorphism of degree one $d$ satisfying
$d^N=0$ ($N\in \mathbb N\backslash \{0,1\}$). In this
case, the $H^{(k)}$ are also $\mathbb Z$-graded and the hexagon
$(\calh^{\ell,m})$ of  Lemma 1 splits into $N$ long exact
sequences $(\cals^{\ell,m}_p),\ p\in\{0,1,\dots,N-1\}$. \[ \begin{array}{ll} &
\cdots\ \longrightarrow  H^{(m),Nr+p}\ \displaystyle{\buildrel {[i^\ell]}\over
\longrightarrow}\ H^{(\ell+m),Nr+p}\ \displaystyle{\buildrel {[d^m]}\over
\longrightarrow}\ H^{(\ell),Nr+p+m\ }\displaystyle{\buildrel
{[i^{N-(\ell+m)}]}\over \longrightarrow}\\  & H^{(N-m),Nr+p+m}\
\displaystyle{\buildrel {[d^\ell]}\over \longrightarrow}\
H^{(N-(\ell+m)),Nr+p+\ell+m}\ \displaystyle{\buildrel {[i^m]}\over
\longrightarrow}\ H^{(N-\ell),Nr+p+\ell+m}\\  & \displaystyle{\buildrel
{[d^{N-(\ell+m)}]}\over \longrightarrow}\ H^{(m),N(r+1)+p}\
\displaystyle{\buildrel {[i^\ell]}\over \longrightarrow}\ \cdots \end{array} \]
where $H^{(k),n}=\{x\in
E^n\vert d^k(x)=0\}/d^{N-k}(E^{n+k-N})$. If instead of being  $\mathbb
Z$-graded, $E$ is $\mathbb Z_N$-graded as in Example 1, then the $N$ exact
sequences
$(\cals^{\ell,m}_p)$ are again $N$ exact hexagons because in $\mathbb Z_N$ one
has $Nr+p=N(r+1)+p=p$.

\section{Graded $q$-differential algebras}

In the remaining part of the paper, $q$ is a complex number with $q\not=0$ and
$q\not=1$.
\begin{definition}
A graded $q$-differential algebra is a $\mathbb N$-graded associative unital
$\mathbb C$-algebra $\fraca=\oplus_n\fraca^n$ equipped with an endomorphism $d$
of degree one satisfying $d(\alpha\beta)=d(\alpha)\beta+q^a\alpha d(\beta),$
$\forall\alpha\in \fraca^a$ and $\beta\in \fraca$, and such that $d^N=0$
whenever $q^N=1$ for $N\in \mathbb N$ with $N\geq 2$.
\end{definition}

The endomorphism $d$ will be referred to as the $q$-{\sl differential of}
$\fraca$ and the above twisted Leibniz rule as the $q$-{\sl Leibniz rule}. For
each $q\in \mathbb C\backslash \{0,1\}$ there is a natural notion of
homomorphism of graded $q$-differential algebra.

An ordinary graded differential algebra is thus a graded $(-1)$-differential
algebra with this terminology. Let us give some other examples.

\subsection*{Example 3: Complex matrix algebras}

Let us come back to Example 1 of Section 2. The only reason why $M_S(\mathbb
C)$ fails to be a graded $q$-differential algebra, with $q$ being a primitive
$N$-th root of the unity, is that it is $\mathbb Z_N$-graded instead of being
$\mathbb N$-graded. There is however an easy way to obtain a graded
$q$-differential algebra from it. Let $p:\mathbb N \rightarrow \mathbb Z_N$ be
the canonical projection and define $(p^\ast M_S(\mathbb C))^n$ for $n\in
\mathbb N$ by $(p^\ast M_S(\mathbb C))^n=(M_S(\mathbb C))^{p(n)}$. Then, on the
$\mathbb N$-graded space $p^\ast M_S(\mathbb C)=\oplus_n (p^\ast M_S(\mathbb
C))^n$ there is a unique product of $\mathbb N$-graded algebra such that the
canonical projection $\pi:p^\ast M_S(\mathbb C)\rightarrow M_S(\mathbb C)$ is
an algebra homomorphism. Furthermore there is a unique endomorphism of degree 1
of $p^\ast M_S(\mathbb C)$, again denoted by $d$, such that $\pi\circ d=d\circ
\pi$. This endomorphism satisfies the $q$-Leibniz rule and $d^N=0$. Thus
$p^\ast M_S(\mathbb C)$ is a graded $q$-differential algebra. $\square$

\subsection*{Example 4: Simplicial forms}

Let $K$ be a simplicial complex, i.e. a set equipped with a set of non-empty
subsets $\fracs$ such that $X\in \fracs$ and $Y\subset X$ implies $Y\in
\fracs,\ (Y\not=\emptyset)$. For $n\in \mathbb N$, an {\sl ordered $n$-simplex}
is a sequence $(x_0,\dots,x_n)$ of elements of $K$ such that
$\{x_0,\dots,x_n\}\in \fracs$. A {\sl simplicial $n$-form} is a complex-valued
function $(x_0,\dots,x_n)\mapsto \omega(x_0,\dots,x_n)$ on the set of ordered
$n$-simplices. Let $\Omega^n_K$ denote the vector space of all simplicial
$n$-forms. The graded vector space $\Omega_K=\oplus_n\Omega^n_K$ is a $\mathbb
N$-graded unital algebra for the product $(\alpha,\beta)\mapsto \alpha\beta$
defined by
$\alpha\beta(x_0,\dots,x_{a+b})=\alpha(x_0,\dots,x_a)\beta(x_a,\dots,x_{a+b})$
for $\alpha\in \Omega^a_K,\  \beta\in\Omega^b_K$ and any ordered
$(a+b)$-simplex $(x_0,\dots,x_{a+b})$. One defines for $q\in \mathbb
C\backslash \{0,1\}$ an endomorphism $d_q$ of degree 1 of $\Omega_K$ by setting
\[
d_q(\omega)(x_0,\dots,x_{n+1})=\sum^n_{k=0}q^k\omega(x_0,\oisk,
x_{n+1})-q^n\omega(x_0,\dots,x_n)
\]
for $\omega\in \Omega^n_K$ and any ordered $(n+1)$-simplex
($x_0,\dots,x_{n+1}$), where $\oiseau$ means omission of $x_k$. Equipped with
$d_q$, $\Omega_K$ is a graded $q$-differential algebra, i.e. $d_q$ satisfies
the $q$-Leibniz rule and $(d_q)^N=0$ whenever $q^N=1$ for $N\in \mathbb N$ with
$N\geq 2$. The above $d_q$ is a particular case of the one of Example 2,
formula (\ref{eq4}), if one takes for the $\fracf_k$ the dual mappings of the
usual simplicial faces. This works as well with $\cala$-valued simplicial forms
where $\cala$ is a unital associative $\mathbb C$-algebra. $\square$

\subsection*{Example 5: Hochschild cochains}

Let $\cala$ be a unital associative $\mathbb C$-algebra. Recall that a
$\cala$-valued Hochschild $n$-cochain is a linear mapping of $\otimes^n\cala$
into $\cala$. The vector space of these $n$-cochains is denoted by
$C^n(\cala,\cala)$. The $\mathbb N$-graded vector space
$C(\cala,\cala)=\oplus_n C^n(\cala,\cala)$ is in a natural way a $\mathbb
N$-graded unital associative algebra with product defined by
$(\alpha\beta)(x_1,\dots,x_{a+b})=\alpha(x_1,\dots,x_a)\beta(x_{a+1},\dots
x_{a+b})$ for $\alpha\in C^a(\cala,\cala),\beta\in C^b(\cala,\cala)$ and
$x_i\in \cala$, where on the right hand side the product is the product of
$\cala$. One defines for $q\in \mathbb C\backslash \{0,1\}$ a linear mapping of
degree 1 of $C(\cala,\cala)$ in itself $\delta_q$ by setting for $\omega\in
C^n(\cala,\cala)$ and $x_i\in \cala$
\begin{equation}
\begin{array}{ll}
\delta_q(\omega)(x_0,\dots,x_n)=x_0\omega(x_1,\dots,x_n)&+
\sum^n_{k=1}q^k\omega(x_0,\dots,x_{k-1}x_k,\dots,x_n)\\
&-q^n\omega(x_0,\dots,x_{n-1})x_n
\end{array}
\label{eq5}
\end{equation}
One has $\delta_q(\alpha\beta)=\delta_q(\alpha)\beta+q^a\alpha\delta_q(\beta)$
for $\alpha\in C^a(\cala,\cala)$ and $\beta\in C(\cala,\cala)$. Furthermore,
$(\delta_q)^N=0$ whenever $q^N=1$ for $N\in \mathbb N$ with $N\geq 2$. In other
words, equipped with $\delta_q$, $C(\cala,\cala)$ is a graded $q$-differential
algebra. Again, $\delta_q$ is a particular case of the $d_q$ of Example 2 in
Section 2, formula (\ref{eq4}), with appropriate cofaces (see in \cite{Lo}).\\
More generally, for an arbitrary bimodule $\calm$ over $\cala$, Formula
(\ref{eq5}) still defines a linear mapping homogeneous of degree 1 on the
graded vector space $C(\cala,\calm)$ of $\calm$-valued Hochschild cochains
which satisfies $\delta^N_q=0$ whenever $q^N=1$ $(N\geq 2)$  and reduces to the
Hochschild coboundary for $q=-1$.\\
When $q$ is a primitive $N$-th root of the unity, the corresponding
$H^{(k),n}$,$(k\in \{1,\dots,N-1\},\ n\in \mathbb N$), are given by
\cite{M.D-V}:  $H^{(k),Nm}=H^{2m},H^{(k),N(m+1)-k}=H^{2(m+1)-1}$, for
$k\in\{1,\dots,N-1\}$ and $m\in \mathbb N$, where the $H^n$ denote the usual
Hochschild cohomology spaces (of $C(\cala,\calm)$), and the other $H^{(k),n}$
vanish. $\square$

\subsection*{Example 6: Dual of a product}

Let $\cala$ be an associative unital $\mathbb C$-algebra and let
$C(\cala)=\oplus_n C^n(\cala)$ be the graded
vector space of multilinear forms on $\cala$; i.e.
$C^n(\cala)=(\otimes^n\cala)^\ast$ is the ($\mathbb C$-) dual of
$\otimes^n\cala$
and $C^0(\cala)=\mathbb C$. By making the natural identifications
$C^n(\cala)\otimes C^m(\cala)\subset C^{n+m}(\cala)$ one sees that $C(\cala)$
is
canonically a $\mathbb N$-graded unital $\mathbb C$-algebra, (the product being
the tensor product over $\mathbb C$). By duality, the product ${\mathfrak
m}:\cala\otimes\cala\rightarrow \cala$ of $\cala$ gives a linear mapping
$\fm^\ast$ of $\cala^\ast$ into $(\cala\otimes \cala)^\ast$ i.e.
$\fm^\ast:C^1(\cala)\rightarrow C^2(\cala)$. For $q\in {\mathbb C}\backslash
\{0,1\}$, $\fm^\ast$ extends into a linear mapping
$\fm^\ast_q:C(\cala)\rightarrow C(\cala)$ which satisfies the graded
$q$-Leibniz
rule with $\fm^\ast_q(C^0(\cala))=0$ and
\begin{equation}
\fm^\ast_q(\omega)(x_0,\dots,x_n)=\sum^n_{k=1}
q^{k-1}\omega (x_0,\dots, x_{k-1} x_k,\dots,x_n) \label{eq6} \end{equation}
for
$\omega\in C^n(\cala)$ with $n\geq 1$ and $x_i\in \cala$. It follows then from
the associativity
of the product of $\cala$ that one has $(\fm^\ast_q)^N=0$ whenever $q^N=1$,
$N\in
{\mathbb N}\backslash \{0,1\}$. Thus $C(\cala)$ equipped with $\fm^\ast_q$ is a
graded $q$-differential algebra.\\
The mapping $\fm^\ast_q$ of (\ref{eq6}) is a particular case of the $\tilde
d_q$ of Example 2 in Section 2, formula (\ref{eq4}), with an obvious choice for
the cofaces and a shift $-2$ in degree.\\
Let $\frach$ be the linear mapping of degree $-1$ of $C^+(\cala)=\oplus_{n\geq
1}C^n(\cala)$ into itself defined by
$\frach(\omega)(x_1,\dots,x_{n-1})=\omega(\bbbone,x_1,\dots,x_{n-1})$ for
$\omega\in C^n(\cala)$ with $n\geq 2$ and by $\frach((C^1(\cala))=0$. Then one
has on $C^+(\cala): \frach\circ \fracm_q^\ast-q\fracm^\ast_q\circ \frach=I$. It
follows then from Lemma 2 that if $q$ is a primitive $N$-th root of the unity
one has $H^{(k),n}=0$ for $k\in \{1,\dots,N-1\}$ and $n\geq 1$. On the other
hand, in this case, one obviously has $H^{(k),0}=\mathbb C\ (k\in
\{1,\dots,N-1\}$). $\square$

\section{Universal $q$-differential calculus}

Our aim in this section is to produce the $q$-analog of the universal
differential envelope of a unital associative $\mathbb C$-algebra $\cala$
\cite{Kar}. We start with a new example of graded $q$-differential algebra
which is itself of interest.

\subsection*{Example 7: The tensor algebra over $\cala$ of $\cala\otimes\cala$}

The tensor
product (over $\mathbb C$) $\cala\otimes\cala$ is in a natural way a bimodule
over $\cala$. The tensor algebra over $\cala$ of the bimodule
$\cala\otimes\cala$
will be denoted by $\mathfrak{T}(\cala)=\oplus_n \mathfrak{T}^n(\cala)$. This
is a unital graded algebra with
$\mathfrak{T}^n(\cala)=\otimes^{n+1}\cala$ and product defined by \[
(x_1\otimes\dots\otimes x_m)(y_1\otimes \dots\otimes
y_n)=x_1\otimes\dots\otimes
x_{m-1}\otimes x_m y_1\otimes y_2\otimes \dots \otimes y_n \ \mbox{for}\
x_i,y_j\in \cala. \]

In particular $\cala$ coincides with the subalgebra $\mathfrak{T}^0(\cala)$.
Since $\cala$ is unital, $\cala\otimes\cala$ is the free bimodule
generated by $\tau=\bbbone \otimes\bbbone$. Hence $\mathfrak{T}(\cala)$ is the
$\mathbb N$-graded algebra generated by $\cala$ in degree 0 and by a free
generator $\tau$ of degree 1. In fact one has $x_0\otimes\dots\otimes
x_n=x_0\tau
x_1\dots \tau x_n$, $\forall x_i\in \cala$.

\begin{lemma} There is a unique linear mapping
$d_q:\mathfrak{T}(\cala)\rightarrow \mathfrak{T}(\cala)$ homogeneous of degree
1
satisfying the $q$-Leibniz rule such that \[ d_q(x)=\bbbone\otimes x-x\otimes
\bbbone=\tau x - x\tau,\ \forall x\in \cala, \] and \[ d_q(\tau)=\tau^2,\
(\mbox{i.e.}\ d_q(\bbbone\otimes \bbbone)=\bbbone \otimes \bbbone \otimes
\bbbone). \] Moreover $d_q$ satisfies $d_q^N=0$ whenever $q^N=1$ for $N\geq 2$,
$N\in \mathbb N$. \end{lemma}

The $q$-differential $d_q$ on $\ft(\cala)$ defined by the above lemma is given
by
\begin{equation}
\begin{array}{ll}
d_q(x_0\otimes\dots\otimes x_n)= &\sum^n_{k=0}q^kx_0\otimes \dots  \otimes
x_{k-1}\otimes \bbbone \otimes x_k\otimes\dots\otimes  x_n\\
& -q^nx_0\otimes\dots\otimes x_n\otimes \bbbone
\label{eq7}
\end{array}
\end{equation}
Thus again this is a particular case of the $d_q$ of Example 2 in Section 2,
formula (\ref{eq4}), with an obvious choice for the cofaces. By induction on
the integer $n$, one obtains the action of the $n$-th power $d^n_q$ of $d_q$ on
$\tau=\bbbone\otimes\bbbone$ and on the elements of $\cala:
d^n_q(\tau)=[n!]_q\tau^{n+1}\ \mbox{and}\ d^n_q(x)=[n!]_q\tau^{n-1}d_q(x),\
\forall x\in \cala.
$

The algebra $\ft(\cala)$ equipped with $d_q$ is a graded $q$-differential
algebra. $\square$

\subsection*{The universal $q$-differential envelope of $\cala$}

Let $\Omega_q(\cala)$ be the smallest $q$-differential subalgebra of
$\ft(\cala)$ equipped with  $d_q$ which contains $\cala$, i.e. the smallest
subalgebra of $\ft(\cala)$, stable by $d_q$, containing $\cala=\ft^0(\cala)$.
This is again a graded $q$-differential algebra which is characterized uniquely
up to an isomorphism by the following universal property, \cite{D-V K},
\cite{M.D-V}.
\begin{theorem} Let $\fraca=\oplus_n\fraca^n$ be
a graded $q$-differential algebra and let $\varphi:\cala\rightarrow \fraca^0$
be
a homomorphism of unital algebras. Then there is a unique homomorphism
$\bar\varphi:\Omega_q(\cala)\rightarrow \fraca$ of graded $q$-differential
algebras inducing $\varphi$. \end{theorem}
It is natural to call $\Omega_q(\cala)$ {\sl the universal $q$-differential
envelope of $\cala$} or {\sl the universal $q$-differential calculus over
$\cala$}. Form $q=-1$, $\Omega_{(-1)}(\cala)$ is just the usual universal
envelope $\Omega(\cala)$ of $\cala$ as defined in \cite{Kar}.\\

Let $q$ be a primitive $N$-th root of the unity and let $E$ be the $\mathbb
Z$-graded vector space $E=\mathbb C e_{-(N-1)}\oplus \dots\oplus \mathbb C
e_{-1}\oplus \ft(\cala)$. One extends $d_q$ to $E$ by setting $d_q
e_{-1}=\bbbone$ and $d_q e_{-k}=e_{-(k-1)}$ for $N-1\geq k\geq 2$. One still
has $d_q^N=0$ on $E$. Let $\omega$ be a linear form on $\cala$ satisfying
$\omega(\bbbone)=1$ and let us define an endomorphism $h$ of degree $-1$ of $E$
by: $h(x_0\otimes\dots\otimes x_n)=\omega(x_0)x_1\otimes\dots\otimes x_n$ for
$x_i\in \cala$ and $n\geq 1$, $h(x_0)=-q^{-1}\omega(x_0)e_{-1}$ for $x_0\in
\cala$, $h(e_{-k})=-q^{-(k+1)}(1+q+\dots+q^k) e_{-(k+1)}$ for $N-2\geq k\geq 1$
and $h(e_{-(N-1)})=0$. One has on $E:h\circ d_q-qd_q\circ h=I$. By Lemma 2,
this implies that $H^{(k)}(E)=0$,  $\forall k\in\{1,\dots,N-1\}$. The subspace
$F=\mathbb C e_{-(N-1)}\oplus\dots\oplus\mathbb C e_{-1}\oplus\Omega_q(\cala)$
is stable by $d_q$ and by $h$, \cite{M.D-V}, therefore $H^{(k)}(F)=0$, $\forall
k\in\{1,\dots,N-1\}$. This implies
$H^{(k),n}(\ft(\cala),d_q)=H^{(k),n}(\Omega_q(\cala))=0$ for $n\geq 1$ and
$H^{(k),0}(\ft(\cala),d_q)=H^{(k),0}(\Omega_q(\cala))=\mathbb C$, $\forall k\in
\{1,\dots,N-1\}$.\\
Thus when $q$ is a primitive $N$-th root of the unity, the corresponding
generalized cohomologies of $(\ft(\cala),d_q)$ and of $\Omega_q(\cala)$ are
trivial. This generalizes a well-known result for $q=-1$.

\end{document}